\documentclass[preprint,12pt]{elsarticle}

\usepackage{amssymb}
\newtheorem{lemma}{Lemma}[section]
\newtheorem{theorem}[lemma]{Theorem}
\newtheorem{corollary}[lemma]{Corollary}

\newtheorem{remark}[lemma]{Remark}

\def\R{\mathbb{R}}

\begin{document}

\begin{frontmatter}



\title{Topological Properties of Strong Solutions for the 3D Navier-Stokes Equations}

\author[label1]{Pavlo~O.~Kasyanov}
\author[label2]{Luisa Toscano}
\author[label1]{Nina~V.~Zadoianchuk}

\address[label1]{Institute for Applied System Analysis,
National Technical University of Ukraine ``Kyiv Polytechnic
Institute'', Peremogy ave., 37, build, 35, 03056, Kyiv, Ukraine,\
kasyanov@i.ua, ninellll@i.ua.}

\address[label2]{University of Naples ``Federico II'', Dep. Math. and
Appl. R.Caccioppoli, via Claudio 21, 80125 Naples,Italy,
luisatoscano@libero.it}

\begin{abstract}
In this paper we give a criterion for the existence of global strong
solutions for the 3D Navier-Stokes system for any regular initial
data.
\end{abstract}


\end{frontmatter}

% \linenumbers
\small
\section{Introduction}

Let $\Omega\subseteq \R^3$ be a bounded open set with sufficiently
smooth boundary $\partial \Omega$ and $0<T<+\infty$. We consider the
incompressible Navier-Stokes equations
\begin{equation}\label{N-STeq:1}
\left\{\begin{array}{l} y_t+(y\cdot \nabla)y=\nu\triangle y -\nabla
p+f \mbox{ in
}Q=\Omega\times(0,T),\\
{\rm div}\ y=0\mbox{ in }Q,\\
y=0\mbox{ on }\partial \Omega\times(0,T),\qquad y(x,0)=y_0(x)\mbox{
in }\Omega,
\end{array}\right.
\end{equation}
where $\nu>0$ is a constant. We define the usual function spaces
\[\mathcal{V}=\{u\in(C_{0}^{\infty}(\Omega))^{3}:{\rm div}\ u=0\},\]
\[
H =\mbox{closure of }\mathcal{V}\mbox{ in }{(L^{2}(\Omega))^{3}},
\quad V=\{u\in (H_{0}^{1}(\Omega))^{3}:\,{\rm div}\ u=0\}.
\]
We denote by $V^*$ the dual space of $V$. The spaces $H$ and $V$ are
separable Hilbert spaces and $V\subset H\subset V^*$ with dense and
compact embedding when $H$ is identified with its dual $H^*$. Let
$(\cdot,\cdot)$, $\|\cdot\|_H$ and $((\cdot,\cdot))$, $\|\cdot\|_V$
be the inner product and the norm in $H$ and $V$, respectively, and
let $\langle\cdot,\cdot\rangle$ be the pairing between $V$ and
$V^*$. For $u,v,w\in V$, the equality
\[
b(u,v,w)=\int\limits_{\Omega}\sum\limits_{i,j=1}^3 u_i\frac{\partial
v_j}{\partial x_i}w_jdx
\]
defines a trilinear continuous form on $V$ with $b(u,v,v)=0$ when
$u\in V$ and $v\in (H_0^1(\Omega))^3$. For $u,v\in V$, let $B(u,v)$
be the element of $V^*$ defined by $\langle
B(u,v),w\rangle=b(u,v,w)$ for all $w\in V$.

We say that the function $y$ is a \textit{weak solution} of Pr.
(\ref{N-STeq:1}) on $[0,T]$, if $y\in L^\infty(0,T;H)\cap
L^2(0,T;V)$, $\frac{dy}{dt}\in L^1(0,T;V^*)$, if
\begin{equation}\label{N-STeq:4}
\frac{d}{dt}(y,v)+\nu ((y,v))+b(y,y,v)=\langle f,v\rangle\quad
\mbox{for all }v\in V,
\end{equation}
in the sense of distributions on $(0,T)$, and if $y$ satisfies the
energy inequality
\begin{equation}\label{N-STeq:5}
V(y)(t)\le V(y)(s)\quad\mbox{for all }t\in [s,T],
\end{equation}
for a.e. $s\in (0,T)$ and for $s=0$, where
\begin{equation}\label{N-STeq:6}
V(y)(t):=\frac12\|y(t)\|_H^2+\nu
\int\limits_{0}^t\|y(\tau)\|_V^2d\tau- \int\limits_{0}^t \langle
f(\tau),y(\tau)\rangle d\tau.
\end{equation}

This class of solutions is called Leray--Hopf or physical one. If
$f\in L^2(0,T;V^*)$, and if $y$ satisfies (\ref{N-STeq:4}), then
$y\in C([0,T];H_w)$, $\frac{dy}{dt}\in L^{\frac{4}{3}}(0,T;V^*)$,
where $H_w$ denotes the space $H$ endowed with the weak topology. In
particular, the initial condition $y(0)=y_0$ makes sense for any
$y_0\in H$.

Let $A:V\to V^*$ be the linear operator associated to the bilinear
form $((u,v))=\langle Au,v\rangle$. Then $A$ is an isomorphism from
$D(A)$ onto $H$ with $D(A)=(H^2(\Omega))^3\cap V.$ We recall that
the embedding $D(A)\subset V$ is dense and continuous. Moreover, we
assume $\|Au\|_H$ as the norm on $D(A)$, which is equivalent to the
one induced by $(H^2(\Omega))^3$. The Problem (\ref{N-STeq:1}) can
be rewritten as
\begin{equation}\label{N-STeqVAr}
\left\{\begin{array}{l} \frac{dy}{dt}+\nu Ay+B(y,y)=f \mbox{ in
}V^*,\\
y(0)=y_0,
\end{array}
\right.
\end{equation}
where the first equation we understand in the sense of distributions on $(0,T)$. Now we write
\[
\mathcal{D}(y_0,f)=\{ y\,: \, y\mbox{ is a weak solution of Pr.
(\ref{N-STeq:1}) on }[0,T]\}.
\]
It is well known (cf. \cite{Temam}) that if $f\in L^2(0,T;V^*)$, and
if $y_0\in H$, then $\mathcal{D}(y_0,f)$ is not empty.

A weak solution $y$ of Pr. (\ref{N-STeq:1}) on $[0,T]$ is called a
\textit{strong} one, if it additionally belongs to Serrin's class
$L^8(0,T;(L^4(\Omega))^3)$. We note that any strong solution $y$ of
Pr. (\ref{N-STeq:1}) on $[0,T]$ belongs to $C([0,T];V)\cap
L^2(0,T;D(A))$ and $\frac{dy}{dt}\in L^2(0,T;H)$ (cf.
\cite[Theorem~1.8.1, p.~296]{Sohr} and references therein).

For any $f\in L^\infty(0,T;H)$ and $y_0\in V$ it is well known the
only local existence of strong solutions for the 3D Navier-Stokes
equations (cf. \cite{Sohr, Temam, ZMK3} and references therein).
Here we provide a criterion for existence of strong solutions for
Pr. (\ref{N-STeq:1}) on $[0,T]$ for any initial data $y_0\in V$ and
$0<T<+\infty$.

\section{Topological Properties of Strong Solutions}

The main result of this note has the following form.

\begin{theorem}\label{mt}
Let $f\in L^2(0,T;H)$ and $y_0\in V$. Then either for any
$\lambda\in [0,1]$ there is an $y_\lambda\in
 C([0,T];V)\cap L^2(0,T;D(A))$ such that
$y_\lambda\in\mathcal{D}(\lambda y_0,\lambda f)$, or the set
\begin{equation}\label{eq:set}
\{ y\in C([0,T];V)\cap L^2(0,T;D(A)) \, : \, y\in\mathcal{D}(\lambda
y_0,\lambda f), \lambda\in(0,1) \}
\end{equation}
is unbounded in $L^8(0,T;(L^4(\Omega))^3)$.
\end{theorem}

In the proof of Theorem~\ref{mt} we use an auxiliary statement
connected with continuity property of strong solutions on parameters
of Pr.~(\ref{N-STeq:1}) in Serrin's class
$L^8(0,T;(L^4(\Omega))^3)$.

\begin{theorem}\label{ml}
Let $f\in L^2(0,T;H)$ and $y_0\in V$. If $y$ is a strong solution
for Pr. (\ref{N-STeq:1}) on $[0,T]$, then there exist $L,\,\delta>
0$ such that for any $z_0\in V$ and $g\in L^2(0,T;H)$, satisfying
the inequality
\begin{equation}\label{N-STeq:7}
\|z_0-y_0\|_V^2+\|g-f\|_{L^2(0,T;H)}^2<\delta,
\end{equation}
the set $\mathcal{D}(z_0,g)$ is one-point set $\{z\}$ which belongs
to $C([0,T];V)\cap L^2(0,T;D(A))$, and
\begin{equation}\label{eq:lip}
\|z-y\|_{C([0,T];V)}^2+\frac{\nu}{2}\|z-y\|_{D(A)}^2\le
L\left(\|z_0-y_0\|_V^2+\|g-f\|_{L^2(0,T;H)}^2\right).
\end{equation}
\end{theorem}
\begin{remark}
We note that from Theorem~\ref{ml} with $z_0\in V$ and $g\in
L^2(0,T;H)$ with $\|z_0\|_V^2+\|g\|_{L^2(0,T;H)}^2$ sufficiently
small, Problem (\ref{N-STeq:1}) has only one global strong solution.
\end{remark}
\begin{remark}
Theorem~\ref{ml} provides that, if for any $\lambda\in [0,1]$ there
is an $y_\lambda\in
 L^8(0,T;(L^4(\Omega))^3)$ such that
$y_\lambda\in\mathcal{D}(\lambda y_0,\lambda f)$, then the set
\[
\{ y\in C([0,T];V)\cap L^2(0,T;D(A)) \, : \, y\in\mathcal{D}(\lambda
y_0,\lambda f), \lambda\in(0,1) \}
\]
is bounded in $L^8(0,T;(L^4(\Omega))^3)$.
\end{remark}

If $\Omega$ is a $C^\infty$-domain and if $f\in
C_0^\infty(\overline{(0,T)\times \Omega})^3$, then any strong
solution $y$ of Pr. (\ref{N-STeq:1}) on $[0,T]$ belongs to
$C^\infty((0,T]\times \Omega)^3$ and $p\in C^\infty((0,T]\times
\Omega)$ (cf. \cite[Theorem~1.8.2, p.~300]{Sohr} and references
therein). This fact directly provides the next corollary of
Theorems~\ref{mt} and \ref{ml}.
\begin{corollary}\label{cor1}
Let $\Omega$ be a $C^\infty$-domain, $f\in
C_0^\infty(\overline{(0,T)\times \Omega})^3$. Then either for any
$y_0\in V$ there is a strong solution of Pr. (\ref{N-STeq:1}) on
$[0,T]$, or the set
\[
\{ y\in C^\infty((0,T]\times \Omega)^3 \, : \,
y\in\mathcal{D}(\lambda y_0,\lambda f), \lambda\in(0,1) \}
\]
is unbounded in $L^8(0,T;(L^4(\Omega))^3)$ for some $y_0\in
C_0^\infty(\Omega)^3$.
\end{corollary}

\section{Proof of Theorem~\ref{ml}} Let $f\in L^2(0,T;H)$,
$y_0\in V$, and $y\in C([0,T];V)\cap L^2(0,T;D(A))$ be a strong
solution of Pr.~(\ref{N-STeq:1}) on $[0,T]$. Due to \cite{Serrin},
\cite[Chapter~3]{Temam} the set $\mathcal{D}(y_0,f)=\{y\}$. Let us
now fix $z_0\in V$ and $g\in L^2(0,T;H)$ satisfying (\ref{N-STeq:7})
with
\begin{equation}\label{delta}
\delta=\min\left\{1;\frac{\nu}{4}\right\}e^{-2TC},\
C=\max\left\{\frac{27c^4}{2\nu^{3}}; \frac{7^8c^8}{2^{12}\nu^{7}}\right\}\left(\|y\|^4_{C([0,T];V)}+1\right)^2,
\end{equation}
$c>0$ is a constant from the inequalities (cf. \cite{Sohr, Temam})
\begin{equation}\label{eq:b1}
|b(u,v,w)|\le
c\|u\|_V\|v\|_{V}^{\frac12}\|v\|_{D(A)}^{\frac12}\|w\|_H \quad
\forall u\in V,\, v\in D(A),\, w\in H;
\end{equation}
\begin{equation}\label{eq:b2}
|b(u,v,w)|\le
c\|u\|_{D(A)}^{\frac34}\|u\|_{V}^{\frac14}\|v\|_{V}\|w\|_H
\quad
\forall u\in D(A),\, v\in V,\, w\in H.
\end{equation}
The auxiliary problem
\begin{equation}\label{N-STeq:8}
\left\{\begin{array}{l} \frac{d\eta}{dt}+\nu
A\eta+B(\eta,\eta)+B(y,\eta)+B(\eta,y)=g-f\mbox{ in
}V^*,\\
\eta(0)=z_0-y_0,
\end{array}
\right.
\end{equation}
has a strong solution $\eta\in C([0,T];V)\cap L^2(0,T;D(A))$ with
$\frac{d\eta}{dt}\in L^2(0,T;H)$, i.e.
\[
\frac{d}{dt}(\eta,v)+\nu
((\eta,v))+b(\eta,\eta,v)+b(y,\eta,v)+b(\eta,y,v)=\langle g-f,v\rangle\quad
\mbox{for all }v\in V,
\]
in the sense of distributions on $(0,T)$. In fact, let
$\{w_j\}_{j\ge 1}\subset D(A)$ be a special basis (cf.
{\cite[p.~56]{Temam1})}, i.e. $Aw_j=\lambda_jw_j$, $j=1,2,...$,
$0<\lambda_1\le \lambda_2\le ... ,\,\, \lambda_j\to +\infty$, $j\to
+\infty$. We consider Galerkin approximations $\eta_m:[0,T]\to
\mbox{span}\{w_j\}_{j=1}^m$ for solutions of Pr. (\ref{N-STeq:8})
satisfying
\[
\frac{d}{dt}(\eta_m,w_j)+\nu
((\eta_m,w_j))+b(\eta_m,\eta_m,w_j)+b(y,\eta_m,w_j)+b(\eta_m,y,w_j)=\langle
g-f,w_j\rangle,
\]
with $(\eta_m(0),w_j)=(z_0-y_0,w_j)$, $j=\overline{1,m}$. Due to
(\ref{eq:b1}), (\ref{eq:b2}) and Young's\,inequality we get
\[
2\langle
g-f,A\eta_m\rangle\le 2\|g-f\|_H\|\eta_m\|_{D(A)}\le \frac{\nu}{4}\|\eta_m\|_{D(A)}^2+\frac{4}{\nu}\|f-g\|_H^2;
\]
\[
-2b(\eta_m,\eta_m,A\eta_m)\le 2c\|\eta_m\|_V^{\frac32}\|\eta_m\|_{D(A)}^{\frac32}\le
\frac{\nu}{2}\|\eta_m\|_{D(A)}^2+ \frac{27c^4}{2\nu^3}\|\eta_m\|_V^6;
\]
\[
-2b(y,\eta_m,A\eta_m)\le 2 c\|y\|_V\|\eta_m\|_V^{\frac12}\|\eta_m\|_{D(A)}^{\frac32}\le 
\frac{\nu}{2}\|\eta_m\|_{D(A)}^2+\frac{27c^4}{2\nu^3}\|y\|_{C([0,T];V)}^4\|\eta_m\|_V^2;
\]
\[
-2b(\eta_m,y,A\eta_m)\le 2c \|\eta_m\|_{D(A)}^{\frac74}\|\eta_m\|_V^{\frac14}\|y\|_{V}\le
\frac{\nu}{2}\|\eta_m\|_{D(A)}^2+\frac{7^8c^8}{2^{12}\nu^7}\|y\|_{C([0,T];V)}^8\|\eta_m\|_V^2.
\]
Thus, 
\[
\frac{d}{dt}\|\eta_m\|_V^2+\frac{\nu}{4}\|\eta_m\|_{D(A)}^2\le
C(\|\eta_m\|_V^2+\|\eta_m\|_V^6)+\frac{4}{\nu}\|g-f\|_H^2,
\]
where $C>0$ is a constant from (\ref{delta}). Hence, the absolutely continuous function
$\varphi=\min\{\|\eta_m\|_V^2,1\}$ satisfies the inequality
$\frac{d}{dt}\varphi\le 2C \varphi+\frac{4}{\nu}\|g-f\|_H^2$, and
therefore $\varphi\le L(\|z_0-y_0\|_V^2+\|g-f\|_{L^2(0,T;H)}^2)< 1$
on $[0,T]$, where $L=\delta^{-1}$. Thus, $\{\eta_n\}_{n\ge 1}$ is
bounded in $L^\infty(0,T;V)\cap L^2(0,T;D(A))$ and
$\{\frac{d}{dt}\eta_n\}_{n\ge 1}$ is bounded in $L^2(0,T;H)$. In a
standard way we get that the limit function $\eta$ of $\eta_n$,
$n\to+\infty$, is a strong solution of Pr. (\ref{N-STeq:8}) on
$[0,T]$. Due to \cite{Serrin}, \cite[Chapter~3]{Temam} the set
$\mathcal{D}(z_0,g)$ is one-point $z=y+\eta \in
L^8(0,T;(L^4(\Omega))^3)$. So, $z$ is strong solution of Pr.
(\ref{N-STeq:1}) on $[0,T]$ satisfying (\ref{eq:lip}).

The theorem is proved.

\section{Proof of Theorem~\ref{mt}}

We provide the proof of Theorem~\ref{mt}. Let $f\in L^2(0,T;H)$ and
$y_0\in V$. We consider the 3D controlled Navier-Stokes system (cf.
\cite{KKV, MT})
\begin{equation}
\left\{
\begin{array}
[c]{l}%
\frac{dy}{dt}+\nu Ay+B(z,y)=f,\\
y(0)=y_{0},
\end{array}
\right.  \label{eq:1}%
\end{equation}
where $z\in L^8(0,T;(L^4(\Omega))^3)$.

By using standard Galerkin approximations (see \cite{Temam}) it is
easy to show that for any $z\in L^8(0,T;(L^4(\Omega))^3)$ there
exists an unique weak solution $y\in L^\infty(0,T;H)\cap L^2(0,T;V)$
of Pr. (\ref{eq:1}) on $[0,T]$, that is,
\begin{equation}
\frac{d}{dt}\left(  y,v\right)  +\nu((y,v))+b(z,y,v)=\left\langle
f,v\right\rangle \mbox{, for all }v\in V, \label{Solution}%
\end{equation}
in the sense of distributions on $(0,T)$.
Moreover, by the inequality%
\begin{equation}\label{eq:b}
|b(u,v,Av)|\le
c_1\|u\|_{(L^4(\Omega))^3}\|v\|_{V}^{\frac14}\|v\|_{D(A)}^{\frac74}\le
\frac\nu2\|v\|_{D(A)}^2+  c_2\|u\|_{(L^4(\Omega))^3}^8\|v\|_{V}^2,
\end{equation}
for all $u \in (L^4(\Omega))^3 \mbox{ and } v\in D(A)$, where
$c_1,c_2>0$ are some constants that do not depend on $u$, $v$ (cf.
\cite{Temam}), we find that $y\in C([0,T];V)\cap L^2(0,T;D(A))$ and
$B(z,y)\in L^{2}(0,T;H)$, so $\frac{dy}{dt}\in L^{2}(0,T;H)$ as
well. We add that, for any $z\in L^8(0,T;(L^4(\Omega))^3)$ and
corresponding weak solution $y\in C([0,T];V)\cap L^2(0,T;D(A))$ of
(\ref{eq:1}) on $[0,T]$, by using Gronwall inequality, we obtain
\begin{equation}\label{eq:Gr}
\begin{array}{c}
\|y(t)\|_V^2\le \|y_0\|_V^2
e^{2c_2\int\limits_0^t\|z(t)\|_{(L^4(\Omega))^3}^8dt},\qquad \forall
t\in
[0,T];\\
\nu\int\limits_0^T \|y(t)\|_{D(A)}^2dt \le \|y_0\|_V^2\left[1+2c_2
e^{2c_2\int\limits_0^T\|z(t)\|_{(L^4(\Omega))^3}^8dt}\|z\|_{L^8(0,T;(L^4(\Omega))^3)}^8\right].
\end{array}
\end{equation}

Let us consider the operator $F:L^8(0,T;(L^4(\Omega))^3)\to
L^8(0,T;(L^4(\Omega))^3)$, where $F(z)\in C([0,T];V)\cap
L^2(0,T;D(A))$ is the unique weak solution of (\ref{eq:1}) on
$[0,T]$ corresponded to $z\in L^8(0,T;(L^4(\Omega))^3)$.

Let us check that $F$ is a compact transformation of Banach space
$L^8(0,T;(L^4(\Omega))^3)$ into itself (cf. \cite{Cr}). In fact, if
$\{z_n\}_{n\ge 1}$ is a bounded sequence in
$L^8(0,T;(L^4(\Omega))^3)$, then, due to (\ref{eq:b}) and
(\ref{eq:Gr}), the respective weak solutions $y_n$, $n=1,2,...$, of
Pr.~(\ref{eq:1}) on $[0,T]$ are uniformly bounded in $C([0,T];V)\cap
L^2(0,T;D(A))$ and their time derivatives $\frac{dy_n}{dt}$,
$n=1,2,...$, are uniformly bounded in $L^{2}(0,T;H)$. So,
$\{F(z_n)\}_{n\ge 1}$ is a precompact set in
$L^8(0,T;(L^4(\Omega))^3)$. In a standard way we deduce  that
$F:L^8(0,T;(L^4(\Omega))^3)\to L^8(0,T;(L^4(\Omega))^3)$ is
continuous mapping.

Since $F$ is a compact transformation of $L^8(0,T;(L^4(\Omega))^3)$
into itself, Schaefer's Theorem (cf. \cite[p.~133]{Cr} and
references therein) and Theorem~\ref{ml} provide the statement of
Theorem~\ref{mt}. We note that Theorem~\ref{ml} implies that the set
$\{z\in L^8(0,T;(L^4(\Omega))^3)\, : \, z=\lambda F(z),\ \lambda\in
(0,1)\}$ is bounded in $L^8(0,T;(L^4(\Omega))^3)$ iff the set
defined in (\ref{eq:set}) is bounded in $L^8(0,T;(L^4(\Omega))^3)$.

The theorem is proved.

 { {{ }
\end{document}